\documentclass[11pt]{amsart}
\usepackage[all]{xy}
\begin{document}
\title{Pseudo-index of Fano manifolds and smooth blow-ups\\
%VERY PRELIMINARY !! DO NOT DISTRIBUTE
}
\author{Laurent BONAVERO}
\date{September 2003}
\maketitle
\noindent
\def\restriction{\string |}
\newcommand{\pp}{\rm ppcm}
\newcommand{\pg}{\rm pgcd}
\newcommand{\Ker}{\rm Ker}
\newcommand{\C}{{\mathbb C}}
\newcommand{\Q}{{\mathbb Q}}
\newcommand{\GL}{\rm GL}
\newcommand{\SL}{\rm SL}
\newcommand{\diag}{\rm diag}
\newcommand{\N}{{\mathbb N}}
\def\finpreuve
{\hskip 3pt \vrule height6pt width6pt depth 0pt}

\newtheorem{theo}{Theorem}
\newtheorem{prop}{Proposition}
\newtheorem{lemm}{Lemma}
\newtheorem{lemmf}{Lemme fondamental}
\newtheorem{defi}{Definition}
\newtheorem{exo}{Exercice}
\newtheorem{rem}{Remark}
\newtheorem{cor}{Corollary}
\newcommand{\CC}{{\mathbb C}}
\newcommand{\ZZ}{{\mathbb Z}}
\newcommand{\RR}{{\mathbb R}}
\newcommand{\QQ}{{\mathbb Q}}
\newcommand{\FF}{{\mathbb F}}
\newcommand{\PP}{{\mathbb P}}
\newcommand{\NN}{{\mathbb N}}
\newcommand{\codim}{\operatorname{codim}}
\newcommand{\Ho}{\operatorname{Hom}}
\newcommand{\Pic}{\operatorname{Pic}}
\newcommand{\NE}{\operatorname{NE}}
\newcommand{\Nun}{\operatorname{N}}
\newcommand{\card}{\operatorname{card}}
\newcommand{\Hilb}{\operatorname{Hilb}}
\newcommand{\mult}{\operatorname{mult}}
\newcommand{\vol}{\operatorname{vol}}
\newcommand{\divi}{\operatorname{div}}
\newcommand{\pr}{\operatorname{pr}}
\newcommand{\con}{\operatorname{cont}}
\newcommand{\ima}{\operatorname{Im}}
\newcommand{\rk}{\operatorname{rk}}
\newcommand{\Exc}{\operatorname{Exc}}
\newcounter{subsub}[subsection]
\def\thesubsub{\thesubsection .\arabic{subsub}}
\def\subsub#1{\addtocounter{subsub}{1}\par\vspace{3mm}
\noindent{\bf \thesubsub ~ #1 }\par\vspace{2mm}}
\def\coker{\mathop{\rm coker}\nolimits}
\def\pr{\mathop{\rm pr}\nolimits}
\def\im{\mathop{\rm Im}\nolimits}
\def\hfl#1#2{\smash{\mathop{\hbox to 12mm{\rightarrowfill}}
\limits^{\scriptstyle#1}_{\scriptstyle#2}}}
\def\vfl#1#2{\llap{$\scriptstyle #1$}\big\downarrow
\big\uparrow
\rlap{$\scriptstyle #2$}}
\def\diagram#1{\def\normalbaselines{\baselineskip=0pt
\lineskip=10pt\lineskiplimit=1pt}   \matrix{#1}}
\def\limind{\mathop{\oalign{lim\cr
\hidewidth$\longrightarrow$\hidewidth\cr}}}

\long\def\InsertFig#1 #2 #3 #4\EndFig{
\hbox{\hskip #1 mm$\vbox to #2 mm{\vfil\includegraphics{#3}}#4$}}
\long\def\LabelTeX#1 #2 #3\ELTX{\rlap{\kern#1mm\raise#2mm\hbox{#3}}}

{\let\thefootnote\relax
\footnote{%\hskip3em 
\textbf{A.M.S.~classification :} 14J45, 14E30, 14E05. 
}}

\vspace{-1cm}

%\begin{center}Pr\'epublication de l'Institut Fourier n$^0$ 493 (2000) \\
%{\em http://www-fourier.ujf-grenoble.fr/prepublications.html}
%\end{center}

\begin{center}
\begin{minipage}{130mm}
\scriptsize

{\bf Abstract. } Suppose $\pi : X \to Y$ is a smooth 
blow-up along a submanifold $Z$ of $Y$ between complex Fano
manifolds $X$ and $Y$ of pseudo-indices $i_X$ and $i_Y$ respectively
(recall that $i_X$ is defined by $i_X := \min 
\{ -K_X \cdot C \,|\, C \mbox{ is a rational curve of } X\}$).
We prove that $i_X \leq i_Y$ if
$2\dim (Z )< \dim(Y) + i_Y -1$
and show that this result is optimal by classifying the 
``boundary'' cases. As expected, these results are obtained
by studying rational curves on $X$ and $Y$.

\end{minipage}
\end{center}

\section{Statement of the results}

\subsection{Introduction} 
When studying surjective morphisms $f~: X \to Y$ between smooth Fano
manifolds $X$ and $Y$ of the same dimension, one generally observes
that the anti-canonical bundle $-K_Y$ of $Y$ is ``more positive'' 
than the anti-canonical bundle $-K_X$ of $X$, one of 
the most important results in this direction being 
the famous theorem of Lazarsfeld \cite{La83} stating that
if $\PP^n \to Y$ is a surjective morphism from $\PP^n$
to an $n$-dimensional manifold $Y$, then $Y \simeq \PP^n$.

For a Fano manifold $X$ (i.e., a complex manifold 
with ample anti-canonical line bundle $-K_X$), one defines two integers
called the index $r_X$ and the pseudo-index $i_X$ of $X$ by 
$$ r_X := \max \{m \in \NN \,| \, -K_X = mL \mbox{ with } L \in \Pic (X)\} $$
and   
$$ i_X := \min \{ -K_X \cdot C \,|\, C \mbox{ is a rational curve of } X\}.$$
Of course, $i_X$ is a multiple of $r_X$ and many results are known 
for these numbers. Among others, Fano manifolds of dimension $n$
with large index 
(namely bigger than $n-2$) are classified (see \cite{IP99} for a complete
survey on Fano manifolds), the situation 
being much more complicated for the pseudo-index: one knows that 
$i_X \leq n +1$ by Mori theory, equality holding if and only
if $X \simeq \PP ^n$ \cite{CMS00}.

\subsection{The main result} Let us start with an easy remark: 
let $Y$ be a complex manifold of 
dimension $n$, let $Z$ be a connected 
submanifold of $Y$, let $X := B_Z(Y)$ be 
the blow-up of $Y$ with center $Z$ and let $E$ be the exceptional divisor
of $\pi : X \to Y$.
We classically have 
$$H^2(X,\ZZ) \simeq H^2(Y,\ZZ) \oplus \ZZ \cdot E
\mbox{ and } K_X = \pi^*K_Y + (n-\dim(Z)-1)E.$$ Therefore, 
if both $Y$ and $X$ are Fano, $r_X$ is equal to the    
greatest common divisor of $r_Y$ and $n-\dim(Z)-1$,
which implies in particular that $r_X \leq r_Y$
and confirms the philosophy described above.

In this Note, we study the behaviour of the pseudo-index
with respect to smooth blow-ups.
Quite surprisingly, this behaviour depends
on the dimension of the center of the blow-up. 
Our precise results are the following.  

%\medskip

\begin{theo}\label{main} Let $Y$ be a complex manifold of 
dimension $n$, let $Z$ be a connected 
submanifold of $Y$ and let $X := B_Z(Y)$ be 
the blow-up of $Y$ with center $Z$.
Suppose both $Y$ and $X$ are Fano. 
\begin{enumerate} 
\item [(i)] If $2\dim (Z )< n + i_Y -1$, then  
$i_X \leq i_Y$,
\item [(ii)] if $2\dim (Z ) = n + i_Y -1$ and 
$i_Y \geq 2$, then $i_X \leq i_Y$,
\item [(iii)] if $\dim (Z )< n/2$,
then $i_X \leq i_Y$. 
\end{enumerate}
\end{theo}

%\medskip

Of course, (iii) is an obvious consequence of (i)
since $i_Y \geq 1$. This result says that the 
pseudo-index has the ``expected behaviour'' when the center
of the blow-up has small dimension. Remark that 
the case where $\dim (Z) =0$ could be proved
by looking at the classification given in \cite{BCW02}
and the case where $\dim (Z) =1$ is Proposition 3.7
of \cite{BCDD03}.

Let us now give an example where $i_X$ is bigger than
$i_Y$. 
In the following proposition (as in the whole 
paper), we do not follow Grothendieck's
convention: $\PP (V)$ denotes the projective space of
{\it lines} of the vector space $V$.  

%\medskip

\begin{prop} Let $n:=2m$ be an even integer,
let ${\mathcal E}$ be the following rank $m+1$ vector bundle
over $\PP^m$:
$$ {\mathcal E} = {\mathcal O}_{\PP^m}^{\oplus m}\oplus {\mathcal O}_{\PP^m}(1)$$
and let $Y_n$ be the $n$-dimensional manifold $Y_n =\PP({\mathcal E})$.
The trivial rank $m$-subbundle of ${\mathcal E}$ defines a submanifold
$Z_m$ isomorphic to $\PP ^{n/2}$ with normal bundle $N_{Z_m/Y_n}$
isomorphic to ${\mathcal O}_{\PP^m}(-1)^{\oplus m}$. Finally, let
$\pi_n : X_n= B_{Z_m}(Y_n) \to Y_n$ be the blow-up of 
$Y_n$ along $Z_m$.
Then $Y_n$ and $X_n$ are Fano manifolds of dimension $n$ if $n \geq 4$.
Moreover $i_{Y_n} = 1$, $i_{X_4} =1$ and $i_{X_n} = 2$ if $n\geq 6$.
\end{prop}

Therefore, the inequalities of Theorem \ref{main} are optimal:
for any $n=2m \geq 6$, $\pi_n : X_n= B_{Z_m}(Y_n) \to Y_n$
is a blow-up with smooth connected center
between Fano manifolds with $\dim(Z_m) = \dim(X_n) /2$ and
$i_{X_n} > i_{Y_n}$.

\medskip

{\em Proof of Proposition 1.} Since $X_n$ and $Y_n$ are naturally toric
manifolds, it is enough to compute the anti-canonical degree
of invariant (rational) curves. If $d$ is a line contained in 
$Z_m$, then $-K_{Y_n} \cdot d = 1$, which gives $i_{Y_n} = 1$.
The Fano manifold $X_n$ is isomorphic to 
the $\PP^1$-bundle 
$\PP _{\PP^{m-1}\times \PP^m} ({\mathcal O}_{\PP^{m-1}\times \PP^m}
\oplus {\mathcal O}_{\PP^{m-1}\times \PP^m}(1,1))$ 
over $\PP^{m-1}\times \PP^m$
hence $i_{X_n}\leq 2$
(the $\PP^1$-fibers having anti-canonical degree equal to $2$).  
Let $E\simeq \PP^{m-1}\times \PP^m$ be the exceptional divisor
of $\pi_n$: the lines contained in a $\PP^{m-1}\times \{*\}\subset E$
have anti-canonical degree equal to $m-1$, hence
$i_{X_4} =1$ and $i_{X_n} = 2$ if $n\geq 6$.\finpreuve  

\medskip

Remark: for $n\geq 8$, the previous computations show that the rational
curves in $X_n$ of minimal anti-canonical degree are not mapped by
$\pi_n$ to curves of minimal anti-canonical degree in $Y_n$.

\medskip

Let us now discuss in more details the optimality
of Theorem \ref{main} by classifying the ``boundary cases''.

\begin{theo}\label{egalite}
Let $\pi: X\to Y$
be a blow-up with smooth connected center $Z$
between Fano manifolds $X$ and $Y$ of dimension $n$.
If $2\dim (Z ) = n + i_Y -1$ then 
$i_X \leq i_Y$ unless $n\geq 6$ is even,  
$X=X_{n}$, $Y=Y_{n}$ and $\pi = \pi_{n}$.
\end{theo}

\subsection{Some consequences} 
The results above have the following consequences
when the pseudo-index of $Y$ is large or
in low dimensions.

\begin{cor}\label{large} Let $\pi: X\to Y$
be a blow-up with smooth connected center
between Fano manifolds $X$ and $Y$ of dimension $n$.
\begin{enumerate} 
\item [(i)] If $i_Y > n/3 -1$, then $i_X \leq i_Y$.
\item [(ii)] If $i_Y = n/3 -1$, then $i_X \leq i_Y$
unless $n=6$, $X=X_6$, $Y=Y_6$ and $\pi = \pi_6$.
\end{enumerate}
\end{cor}
 
{\em Proof of Corollary \ref{large}.}

{\it Proof of (i).} Suppose by contradiction that 
$i_X > i_Y$. Then by Theorem \ref{main}(i),
$2\dim (Z ) \geq n + i_Y -1$.
But the lines contained in the non-trivial fibers 
of the blow-up are rational curves of
anti-canonical degree $n-1-\dim(Z)$, therefore
$n-1-\dim(Z) \geq i_X > i_Y$, 
hence 
$$n-1-i_Y \geq \dim(Z)+1 \geq n/2 +i_Y /2 +1/2$$
and $i_Y \leq n/3 -1$, a contradiction.

{\it Proof of (ii).} Suppose that $i_Y = n/3 -1$
and that $i_X > i_Y$. 
The previous computations implies that every inequality 
occuring in the proof of (i) is an equality.
In particular, one has $2\dim (Z ) = n + i_Y -1$.
Therefore, Theorem \ref{egalite} implies that $n$ is even
and $Y = Y_{n}$. In particular, $i_Y =1 =  n/3 -1$,
hence $n=6$, which ends the proof. \finpreuve

\medskip

Remark that according to the generalised Mukai conjecture,
as stated and studied in \cite{BCDD03}, Fano manifolds 
$Y$ of dimension $n \geq 6$ with $i_Y > n/3 -1$
should have Picard number $\rho_Y$ satisfying
$\displaystyle{ \rho_Y < \frac{3n}{n-6}}$.
Corollary \ref{large} has the immediate 
following corollary.

\begin{cor}\label{dimpetites} Let $\pi: X\to Y$
be a blow-up with smooth connected center
between Fano manifolds $X$ and $Y$ of dimension $n$.
\begin{enumerate} 
\item [(i)] If $n \leq 5$, then $i_X \leq i_Y$. 
\item [(ii)] If $n =6$, then $i_X \leq i_Y$ unless 
$X=X_6$, $Y=Y_6$ and $\pi = \pi_6$.

\end{enumerate}
\end{cor}

\section{Proofs}

\subsection{Proof of Theorem \ref{main}} It is enough to prove
assertion (i), since (iii) is an obvious consequence of (i)
and (ii) is an immediate consequence of Theorem \ref{egalite}
(note that the Fano manifolds $Y_n$ have pseudo-index $1$).

Let $\pi: X=B_{Z}(Y) \to Y$
be a blow-up with smooth center $Z$ between
Fano manifolds $X$ and $Y$. We will denote by $E = \pi^{-1}(Z)$
the exceptional divisor of $\pi$.
The basic idea is very simple: we take a rational curve
$C$ in $Y$ such that $-K_Y \cdot C = i_Y$ and we want 
to show that there is a rational curve 
$\tilde C$ in $X$, mapping surjectively to $C$ by $\pi$,
such that $-K_X \cdot \tilde C \leq i_Y$.

Suppose first that there is a rational curve
$C$ in $Y$ such that $-K_Y \cdot C = i_Y$ 
and such that $C$ is not contained in $Z$.
The strict transform $\tilde C$ of $C$ is a rational curve
satisfying $E\cdot \tilde C \geq 0$ and 
the formula $-K_X = \pi^*(-K_Y) - (n-1-\dim(Z))E$
immediately implies that $-K_X \cdot \tilde C \leq i_Y$.

Now take a rational curve
$C$ in $Y$ such that $-K_Y \cdot C = i_Y$ and assume 
$C \subset Z$. Let us decompose $N_{Z/Y\, |C}$ 
(we allow here a slight abuse of notations, since
$C$ might be a singular rational curve, we should rather
write $\nu ^*(N_{Z/Y\, |C})$ where $\nu : \PP^1 \to C$ is the normalisation
of $C$):
$$ N_{Z/Y\,|C} = \bigoplus_{i=1}^r {\mathcal O}_{\PP^1}(a_i)$$ 
where $r = n-\dim(Z)$.
The sub-line bundles ${\mathcal O}_{\PP^1}(a_i)$
of $N_{Y/Z \,|C}$ define rational curves $\tilde C_i$ of $E$
satisfying 
$ -K_X \cdot \tilde C_i = -K_Y\cdot C - (r-1)a_i$. 
Therefore, we are done if there exists $i$ such that $a_i \geq 0$.
Suppose the contrary, namely that $a_i \leq -1$ for all $i$.  
Then 
$$-K_Z \cdot C = -K_Y \cdot C - \deg (N_{Z/Y \,|C})
\geq i_Y + \rk (N_{Z/Y \,|C}) = i_Y + n - \dim(Z).$$
Under the assumption $2\dim (Z )< n + i_Y -1$,
we get 
$ -K_Z \cdot C > \dim (Z) +1$ hence, by Mori's 
bend-and-break lemma (see for example \cite{Deb01}, p. 58),
the curve $C$ is numerically equivalent in $Z$, hence in 
$Y$, to a connected
nonintegral effective rational $1$-cycle (passing through 
$2$ arbitrary fixed points of $C$). Each reduced irreducible
component of this $1$-cycle has $-K_Y$ anti-canonical degree
strictly less than $i_Y$, contradiction!~\finpreuve 

\medskip

Remarks: in the above proof, we only used that
$-K_Y \cdot C >0$ for any rational curve of $Y$.
Moreover, the previous proof also shows the following. 
Let $Y$ be a complex manifold, let $Z$ be a connected 
submanifold of $Y$ and let $X := B_Z(Y)$ be 
the blow-up of $Y$ with center $Z$.
Suppose both $Y$ and $X$ are Fano and $i_X > i_Y$.
Then any rational curve $C$ satisfying $-K_Y \cdot C = i_Y$
is contained in $Z$ and for any such curve $C$,
the vector bundle $N_{Z/Y\,|C}^*$ is ample.

\subsection{Proof of Theorem \ref{egalite}} This proof 
assumes that the reader has some familiarity with Mori theory,
see for example \cite{Deb01} for a nice introduction.

By the remark above,
if $i_X > i_Y$, any rational curve
$C$ in $Y$ such that $-K_Y \cdot C = i_Y$ is contained
in $Z$. For such a curve $C$, which has minimal degree with respect to 
an ample line bundle, its deformations in $Z$ containing a given point 
cover a subvariety of dimension $\geq -K_Z\cdot C - 1$  
(recall that this an easy consequence of Riemann-Roch   
formula and the bend-and-break lemma, see for example
\cite{Deb01}, \S 6.5).
Since the computations 
in the proof of Theorem \ref{main}(i)
show that $-K_Z \cdot C = \dim(Z) +1$, 
the deformations of $C$ in $Z$
containing a given point
cover $Z$. Therefore, the Picard number of
$Z$ is one (see \cite{Ko96}, IV 3.13.3). One deduces that
any rational curve $C'$ of $Z$ 
is numerically proportional (in $\Nun _1(Z)$) to $C$, and since
$C$ has minimal anti-canonical degree in $Y$, $C'$
satisfies $-K_Z \cdot C' \geq \dim(Z) +1$, hence
$Z \simeq \PP^{\dim(Z)}$ by \cite{CMS00}. Finally, 
the computations above also show that
for any line $d$ in $Z$,
$N_{Z/Y \, |d} \simeq {\mathcal O}_{d}(-1)^{\oplus n-\dim(Z)}$,
which implies 
that 
$N_{Z/Y} \simeq {\mathcal O}_{\PP^{\dim(Z)}}(-1)^{\oplus n-\dim(Z)}$
by Theorem (3.2.1) in \cite{OSS81}.

Let $E = \PP^{\dim(Z)} 
\times \PP^{n-\dim(Z)-1}$ be the exceptional divisor of $\pi$,
and let $\omega$ be a Mori extremal rational curve in $X$ such that
$E \cdot \omega >0$ (such a curve exists by the classical following
argument: take any curve with strictly positive intersection 
with $E$ and decompose it in the Mori cone $\NE(X)$ 
as an effective combination of extremal curves, at least one of these
curves has strictly positive intersection 
with $E$). 
The corresponding Mori contraction
$\varphi _{\omega}$ satisfies Wi\'sniewski's inequality \cite{Wi91}:
$$ \dim(\Exc(\varphi _{\omega})) + \dim (f)
\geq n -1 + i_X$$
where $\Exc(\varphi _{\omega})$ is the locus of contracted curves,
and $f$ is any non-trivial fiber of $\varphi _{\omega}$.
Since every contracted curve is proportional to $\omega$ in $\Nun _1(X)$
and since ${\mathcal O}(E)_{|E} \simeq 
{\mathcal O}_{\PP^{\dim(Z)} 
\times \PP^{n-\dim(Z)-1}}(-1,-1)$, none
of these curves are contained in $E$, therefore any
non-trivial fiber $f$ of $\varphi _{\omega}$ satisfies $\dim(f) =1$.
Moreover, $i_X \geq 2$, therefore $\Exc(\varphi _{\omega})= X$ 
by Wi\'sniewski's inequality above and
$\varphi _{\omega}$ is 
a fibration which, by Ando's classification \cite{An85}, is a smooth 
$\PP^1$-bundle over an $(n-1)$-dimensional 
Fano manifold $X'$. Moreover, $i_X=2$, hence $i_Y=1$,
therefore $\dim(Z) = n/2$. 

Let us show now that $E\cdot f =1$ for any fiber of $\varphi _{\omega}$.
Indeed, if $\PP^1 \to X'$ is a rational curve of $X'$,
the surface $S = \PP^1 \times_{X'} X$ is a ruled surface,
i.e., a Hirzebruch surface, and the exceptional curve
of $S$ is nothing else than $\PP^1 \times_{X'} E$, which is a section
of $S \to \PP^1$. Hence $E\cdot f =1$ for any fiber of $\varphi _{\omega}$.

One immediately deduces that 
$\varphi _{\omega}: E \to X'$ is an isomorphism,
hence $X = \PP({\mathcal E})$ for some rank $2$ bundle ${\mathcal E}$ over
$X' \simeq E \simeq \PP^{n/2}\times \PP^{n/2 -1}$, and $E$ defines 
a sub-line bundle of ${\mathcal E}$. Therefore ${\mathcal E}$
splits and since 
$N_{E/X} \simeq {\mathcal O}_{\PP^{n/2}\times \PP^{n/2 -1}}
(-1,-1)$, one deduces that
$$ X \simeq \PP({\mathcal O}_{\PP^{n/2}\times \PP^{n/2 -1}}\oplus  
{\mathcal O}_{\PP^{n/2}\times \PP^{n/2 -1}}(1,1)),$$
which ends the proof. \finpreuve

\section{Some comments and some more examples}

\subsection{On the normal bundle of the center.}

The following proposition sheds some light on the example
explained in Proposition 1.

\begin{prop}
Let $\pi: X\to Y$
be a blow-up with smooth connected center $Z$
between Fano manifolds $X$ and $Y$ of dimension $n$.
Suppose moreover that the conormal bundle $N_{Z/Y}^*$
is ample.
Then, 
\begin{enumerate} 
\item [(i)] either $i_X \leq i_Y$,
\item [(ii)] or $i_X =2$, $i_Y=1$, $Z$ is a Fano manifold, 
$X$ is a $\PP^1$-bundle over the $(n-1)$-dimensional Fano manifold
$\PP (N_{Z/Y})$ and $Y$ is a $\PP^{n-\dim(Z)}$-bundle over 
$Z$.
\end{enumerate}
\end{prop}

{\em (Sketch of) proof.} Suppose $i_X > i_Y$ and 
let denote by $E$ the exceptional divisor of $\pi$. 
Since $N_{Z/Y}^*$ is ample, $-E_{|E}$ is also ample
and by Grauert's criterion, $E$ is contractible to a point.
Moreover, if $\omega$ is a Mori extremal rational curve in $X$ such that
$E \cdot \omega >0$, using the same arguments as in the proof
of Theorem \ref{egalite},
the corresponding Mori contraction
$\varphi _{\omega}$ is a $\PP ^1$-bundle over the Fano 
manifold $E \simeq \PP (N_{Z/Y})$.
Finally,
$Z$ is Fano by \cite{SW90} and
one has $\rho_Y +1 = \rho_X = \rho_E +1 = \rho_Z +2$,
hence $\rho_Y = \rho_Z +1$.
This implies that there is at least one Mori extremal curve
of $Y$ which is not contained in $Z$. Since $\pi$ is surjective,
the Mori cone $\NE (Y)$ is generated by the images of 
Mori extremal curve of $X$, which are contained in $E$, except
for the fibers $f$ of the $\PP ^1$-bundle structure $X \to E$.
This implies that $\pi(f)$ is extremal in $Y$ and the corresponding
extremal contraction $\psi: Y \to W$ is a fibration. 
But then, the fibers of $\psi$ have dimension less or equal 
to $n - \dim(Z)$, hence equal to $n - \dim(Z)$ since
$-K_Y \cdot \pi(f) = n - \dim(Z)+1$. Thereore the generic fiber of $\psi$
is $\PP ^{n-\dim(Z)}$, and finally every fiber of $\psi$ 
is $\PP ^{n-\dim(Z)}$ and meets $Z$ transversally at exactly one point
(all this is verified since on $X$, the extremal contraction
associated to $f$ is a $\PP^1$-bundle).
Finally, $W \simeq Z$.\finpreuve

\subsection{Some examples.}\label{ex}

The previous Proposition implies that
if $\pi: X\to Y$ 
is a blow-up
with smooth connected center $Z$
between Fano manifolds $X$ and $Y$ 
with $i_X > i_Y \geq 2$,
then the conormal bundle $N_{Z/Y}^*$ is not ample,
althought its restriction to any rational curve of
minimal $-K_Y$-degree (recall that such a curve 
has to be contained in $Z$) is ample as we saw in the proof
of Theorem \ref{main}~!
It is therefore the good place to give a list of examples
(communicated to me by Cinzia Casagrande) 
of blow-ups $\pi: X\to Y$
with smooth connected center $Z$
between Fano manifolds $X$ and $Y$
with $i_X > i_Y \geq 2$.

\medskip

{\bf Examples.} Let $a$, $d$, $r$ and $s$ be positive integers,
let ${\mathcal E}$ be the following rank $r+s$ vector bundle
over $\PP^a$:
$ {\mathcal E} = 
{\mathcal O}_{\PP^a} ^{\oplus r} \oplus 
{\mathcal O}_{\PP^a}(d)^{\oplus s}$
and let ${\mathcal F} = {\mathcal O}_{\PP^a}(d)^{\oplus s}$ be the rank $s$
vector subbundle of ${\mathcal E}$ defined by the 
${\mathcal O}_{\PP^a}(d)$'s factors.
Define $Y := \PP ({\mathcal E})$ and $Z:= \PP({\mathcal F}) 
\simeq \PP^a\times \PP^{s-1}$. This submanifold $Z$ of $Y$ has
codimension $r$ and normal bundle
$N_{Z/Y}$ equal to ${\mathcal O}_{\PP^a\times \PP^{s-1}}(-d,1)^{\oplus r}$.
Finally, let $X:= B_Z(Y)$ be the blow-up of $Y$ with center $Z$.
Easy computations show that $Y$ is Fano if and only if $a\geq rd$
and $X$ is Fano if and
only if $a \geq d$, and that under these assumptions,
one has 
$$ i_Y = \min(r+s, 1+a-rd) \mbox{ and } i_X = \min(r-1, s+1, 1+a-d),$$ 
which leads 
to many examples satisfying $i_X > i_Y \geq 2$. 
The example of lowest dimension
(namely $10$)
for such $X$ and $Y$ is given when $(a,d,r,s) = (5,1,4,2)$.
Let us also say that many of these examples 
lead to $X$ and $Y$ satisfying $r_X < i_X$ and $r_Y < i_Y$.
In the case $s=1$, the Fano manifolds $Y$'s have been 
considered by Debarre (see \cite{Deb01}, \S 5.11) to construct Fano manifolds
of high degree $(-K_Y)^{\dim (Y)}$.

\subsection{Minimal degree of free rational curves.} 
When studying Fano manifolds, one often 
uses {\em free} rational curves, which means 
rational curves $f: \PP^1 \to X$ such that 
$$f^* T_X \simeq {\mathcal O}_{\PP^1}(a_1)\oplus \cdots \oplus 
{\mathcal O}_{\PP^1}(a_{\dim X})$$ with all the $a_i$'s greater
or equal to $0$ (see \cite{Deb01}, Chapter 4 for details). 
One may then introduce another invariant: 
$$f_X := 
\min \{ -K_X \cdot C \,|\, C \mbox{ is a free rational curve of } X\},$$
which is of great importance in Hwang and Mok's recent works.
It can be interpreted as the minimal anti-canonical degree of rational
curves whose deformations cover an open dense subset of $X$.
It is an easy exercise to show that 
if $f:X \to Y$ is a surjective morphism between Fano
manifolds $X$ and $Y$, then $f_X \leq f_Y$. 
Of course, in any of the examples above where $i_X > i_Y$,
the rational curves in $Y$ of minimal anti-canonical degree are not
free curves and their deformations do not cover any dense open subset
of $Y$. 
  
\subsection{A final remark on a related question.}
In the above results, the assumption that both $X$ and $Y$
are Fano is essential:
when $\pi: X\to Y$
is a blow-up with smooth connected center $Z$
between complex manifolds $X$ and $Y$, understanding 
on which conditions $X$ Fano (resp. $Y$ Fano) implies $Y$ Fano
(resp. $X$ Fano) is a completely different question, whose study 
has been initiaded by Wi\'sniewski in \cite{Wi91}. In particular, no condition on 
the dimension of the center is neither necessary nor sufficient 
(except of course when $Z$ is a point, see \cite{BCW02} for
a complete classification) to get one of the implications above:
the examples in \S \ref{ex}, in the particular case where $rd > a \geq d$, 
give examples              
of smooth blow-up $\pi: X\to Y$ 
between complex manifolds $X$ and $Y$ with $X$ being Fano and $Y$ not.

\medskip

{\em Thanks to Cinzia Casagrande and Olivier Debarre for their comments
on a preliminary version of this Note.}

-----------

{\em 
\noindent Laurent Bonavero. Institut Fourier, UMR 5582, 
Universit\'e de Grenoble 1, BP 74. 
38402 Saint Martin d'H\`eres. FRANCE\\
\noindent e-mail : bonavero@ujf-grenoble.fr
}

\end{document}